 \documentclass[11pt]{amsart}

\usepackage{amssymb}
\usepackage{lineno,hyperref} 
\usepackage{graphicx}
\usepackage{mathrsfs}
\usepackage{amsmath,amsthm,amssymb,enumerate}
\usepackage{xparse,aliascnt,hyperref,bookmark}
\hypersetup{
     colorlinks = true,
     linkcolor = blue,
     anchorcolor = blue,
     citecolor = blue,
     filecolor = blue,
     urlcolor = blue
     }

\usepackage{epsfig}  		
\usepackage{epic,eepic}       
\NewDocumentCommand{\xnewtheorem}{m o m}
 {%
  \IfNoValueTF{#2}
   {\newtheorem{#1}{#3}}
   {%
    \newaliascnt{#1}{#2}%
    \newtheorem{#1}[#1]{#3}%
    \aliascntresetthe{#1}%
    \expandafter\newcommand\csname #1autorefname\endcsname{#3}%
   }%
 }
 \usepackage{enumerate}
 \usepackage{color}

\xnewtheorem{theorem}{Theorem}[section]
\xnewtheorem{proposition}[theorem]{Proposition}
\xnewtheorem{lemma}[theorem]{Lemma}
\xnewtheorem{corollary}[theorem]{Corollary}

\theoremstyle{definition}
\xnewtheorem{definition}[theorem]{Definition}
\xnewtheorem{example}[theorem]{Example}

\xnewtheorem{note}[theorem]{Note}

\theoremstyle{remark}

\xnewtheorem{remark}[theorem]{Remark}

\numberwithin{equation}{section}

\newcommand{\R}{\mathbb{R}}  

\newcommand{\C}{\mathbb{C}}
\newcommand{\N}{\mathbb{N}}

\begin{document}


\title[Linear operators preserving entire functions]{Linear $q$-difference, difference and differential operators
	preserving some $\mathcal{A}$-entire functions}


\author{Jiaxing Huang}
\address{College of Mathematics and Statistics, Shenzhen University, Shenzhen, PR China}
\email{hjxmath@163.com}


\author{Tuen-Wai Ng}
\address{Department of Mathematics, The University of Hong Kong, 
Pokfulam Road, Hong Kong}
\email{ntw@maths.hku.hk}



\begin{abstract}
	We apply Rossi's half-plane version of Borel's Theorem to study the zero distribution of linear combinations of  $\mathcal{A}$-entire functions
(Theorem \ref{thm:FD0}). This provides a unified way to study linear $q$-difference, difference and differential operators (with entire coefficients)
preserving subsets of $\mathcal{A}$-entire functions, and hence obtain several analogous results for the Hermite-Poulain Theorem to linear finite ($q$-)difference operators with polynomial coefficients. The method also produces a result on the existence of infinitely many non-real zeros of some differential polynomials of functions in certain sub-classes of $\mathcal{A}$-entire functions.


\end{abstract}



 \maketitle



\section{Introduction and main results}\label{sec1}
The investigations of linear operators preserving real-rootedness of certain classes of entire functions of one complex variable has a long history.  
In the 1870s, the linear operator preserving the class of {\it hyperbolic polynomials} $\mathcal{HP}$ (i.e. polynomials with real coefficients whose zeros are all real) was initiated by Hermite, and further developed by Laguerre.
In 1914, P\'olya and Schur \cite{PS14} completely described the operators acting diagonally on the standard monomial basis $1, x, x^2, \dots, $ of $\R[x]$ and preserving $\mathcal{HP}$. 
One may then consider the corresponding classification problem for some classes of entire functions containing $\mathcal{HP}$, 
for example, the classical {\it Laguerre-P\'olya class} (see \cite[Definition 5.4.11]{RS}).

Let $S$ be a subset in the complex plane $\C$. An entire function $f$ is said to be in the \emph{S-Laguerre-P\'olya class}, $f\in \mathcal{LP}(S)$, if
\begin{equation}\label{eqn:f}
	f(z)=h(z)e^{-\alpha z^2+ \beta z}, \quad h(z)=cz^n\prod_{k=1}^{\infty}(1-z/z_k)e^{t z/z_k}
\end{equation} 
where $\beta\in\R$, $c\in\R\setminus \{0\}$, $\alpha\geq 0$, $t=\{0, 1\}$, 
$n$ is a non-negative integer and $\{z_k\}$ is a finite or infinite sequence in $S$ with $\sum_k |z_k|^{-t-1}<\infty$. 
By \cite[Theorem 1.11]{Hay64} or \cite[Theorem 3.8.5]{Lan07}, 
for $M(r, h)=\displaystyle\max_{|z|=r}|h(z)|$, we have  $\log M(r, h)=o(r^{t+1})$ as $r\rightarrow\infty$.

Any function $f$ in 
\emph{S-Laguerre-P\'olya class} with $t=0$, $\alpha=0$ and $\beta\geq 0$ is said to be of $S$-type I. 
Clearly, \emph{S-Laguerre-P\'olya class} and its sub-class $S$-type I are generalisations of the classical Laguerre-P\'olya class (when $S=\R$ and $t=1$) 
and type I class (when $S=\R_-:=\{x\in\R| x<0\}$ or $S=\R_+:=\{x\in\R| x>0\}$) respectively.
Notice that by a theorem of P\'olya \cite{Pol13} (see \cite[Theorem 5.4.12]{RS}), $\mathcal{LP}(\R)$ is the closure (in the sense of the uniform convergence on compacta) of polynomials in   $\mathcal{HP}$.

To describe our results, we need to introduce some basic definitions in the Nevanlinna theory (see \cite{Hay64}) and the $\mathcal{A}$-entire functions (which are called class $\mathcal{A}$ functions in Chapter V of B. Ja. Levin's book \cite{Lev64}).

\begin{definition}
	A sequence $\{a_n\}$ of $\C$ is called an {\it $A$-sequence} if it satisfies
	the condition
	\begin{equation}\label{eqn:A}
		\sum_{n=1}^{\infty}\left|\mathrm{Im} \frac{1}{a_n}\right|<\infty.
	\end{equation}
	Here, if $a_n=0$, then we define $\mathrm{Im}\dfrac{1}{a_n}=0$. 
	An entire function $f$ is in class $\mathcal{A}$ if its zero set $\{z_n\}$ is an $A$-sequence and we will also call such $f$ an $\mathcal{A}$-entire function.
	If such $f$ is also in $\mathcal{LP}(\{z_n\})$ or $\{z_n\}$-type I class, then we say $f$ is in class $\mathcal{LP}(A)$ or $A$-type I class respectively.
	Finally, $\mathcal{LP}(A;2)$  is the subset of $\mathcal{LP}(A)$ which contains $f$ with $\alpha >0$ in (\ref{eqn:f}) and $\mathcal{LP}_{t_0}(A;2)$ ($\mathcal{LP}_{t_0}(A)$) is the subset of $\mathcal{LP}(A;2)$ ($\mathcal{LP}(A)$) which contains $f$ with $t=t_0$ in (\ref{eqn:f}).
\end{definition}

Clearly, class $\mathcal{A}$ contains entire functions with only real zeros and hence 
$\mathcal{HP}\subset \mathcal{LP}(\R)\subset\mathcal{A}$. 
Also, by definition, we have $\mathcal{LP}_{t_0}(A;2) \subset \mathcal{LP}(A;2) \subset \mathcal{LP}(A) \subset \mathcal{A}$.  

For any meromorphic function $f$ on $\mathbb{C}$, let $n(r,f)$ be the number of poles of $f$ in $|z|\leq r$, and  $$N(r,f):=\int_1^r\frac{n(t,f)-n(0,f)}{t}dt+n(0,f)\log r$$
be the {\it counting function} of $f$. 
The {\it proximity function} $m(r, f)$ is defined by $$m(r, f):=\frac{1}{2\pi}\int_0^{2\pi}\log^+|f(re^{i\theta})|d\theta,$$ where $\log^+a:=\max(\log a, 0)$.
The {\it Nevanlinna  characteristic function} $T(r,f)$ is defined by $$T(r, f):=m(r,f)+N(r, f).$$ 
We also introduce the {\it exponent of convergence of the zeros} of $f$, $\lambda(f)$, and the {\it order} $\rho(f)$ of $f$, which are given respectively by 
\begin{equation*}\label{eqn:ce}
	\lambda(f):=\limsup_{r\rightarrow\infty}\frac{\log N(r,1/f)}{\log r}=\limsup_{r\rightarrow\infty}\frac{\log n(r,1/f)}{\log r}.
\end{equation*} and 
$$\rho(f):=\limsup_{r\rightarrow\infty}\frac{\log T(r, f)}{\log r}.$$
It is clear that $\lambda(f)\leq \rho(f)$. If $f$ is an entire function, it is not hard to see that  $T(r, f)\leq\log^+ M(r, f)$, for all $r>0$. 
Finally, by $\mathcal{M}_k$ we mean the field  of meromorphic functions 
$f$ with $T(r, f) = o(r^{k})$ as $r \to \infty$ outside a
set of finite measure and by $\C[z]$ we mean the ring of polynomials with complex number coefficients. The field $\mathcal{M}_k$ appears naturally in the studies of some hypertranscendental functions (see for example \cite{HuangNg2020}).

Using some ideas from Eremenko and Rubel \cite{ER96} and Ng and Yang \cite{NgY97}, and the half-plane version of Borel's lemma by Rossi \cite{Ro88} (see Lemma \ref{lem:abor}),  we obtain the following
\begin{theorem}\label{thm:FD0}
	Let $f_1, \ldots, f_n$ be linearly independent entire functions over  $\mathcal{M}_2$ satisfying $N(r, 1/f_i)=o(r^2)$ for all $i$.  Let $a_1,\ldots, a_n$ be entire functions in $\mathcal{M}_1$. Suppose that each $f_i$ is in class $\mathcal{A}$. Then 
	$$F=a_1f_1+\dots+a_nf_n$$ is in $\mathcal{A}\backslash\{0\}$
	if and only if $a_i\not\equiv 0$ for at most one $i$ and this $a_i$ is in class $\mathcal{A}\backslash\{0\}$.
\end{theorem}

\subsection{Linear operators preserving real-rootedness}\label{Sec:L}

We now explain how Theorem \ref{thm:FD0} can be used to classify linear ($q$-difference, difference, differential) operators of finite order preserving some sub-classes of $\mathcal{A}$-entire functions.

Let $\mathcal{M}$ be the field of meromorphic functions on $\C$. Consider the linear difference operator $\Delta_{M_1, M_2, h}: \mathcal{M}\rightarrow\mathcal{M}$ defined by
\begin{equation}\label{eqn:M12}
	\Delta_{M_1, M_2, h}(f)(z):=M_1(z)f(z+ih)+M_2(z)f(z-ih),
\end{equation}
where $M_1$ and $M_2$ are complex-valued functions, and $h$ is a complex number. 
In 1926,  P\'olya \cite{Pol26} established that $\Delta_{1, 1,c}(\mathcal{LP}(\R))\subset\mathcal{LP}(\R)$ for every real number $c$ 
and de Bruijn \cite{Bru50} noticed that actually $\Delta_{e^{i\theta}, e^{-i\theta}, c}(\mathcal{LP}(\R))\subset\mathcal{LP}(\R)$ for every real number $c$ and $\theta \in [0,2\pi]$. 
These studies were continued by Oberschkov \cite{Ob63}, Levin \cite{Lev64}, Craven and Csordas \cite{CC89, CC96}, Walker \cite{Wal97} and others. 
On the other hand, Br\"and\'en and Borcea \cite{BB09} have completely characterised all linear operators on $\C[z]$ preserving $\mathcal{HP}$ as well as $\mathcal{LP}(\R)$ in the entire function space. 
However, it is not easy to apply their characterisation to check when the two term difference operator (\ref{eqn:M12}) preserves $\mathcal{LP}(\R)$. 
In 2017, Katkova et al. \cite{KTV17} thoroughly solved this classification problem and hence generalised the results of P\'olya and de Bruijn. 
They gave the necessary and sufficient conditions of the operator (\ref{eqn:M12}) preserving the class $\mathcal{LP}(\R)$ in terms of the explicit expressions of $M_1$ and $M_2$ (see \cite[Theorem 1.1]{KTV17}).

It is then natural to consider linear difference operators with more than two terms. When the coefficients are polynomials, by applying Br\"and\'en and Borcea's characterisation of linear operators on $\C[z]$ 
preserving $\mathcal{HP}$ (\cite[Theorem 1]{BB09}), Br\"and\'en et al. \cite{BKS16} obtained the following
\begin{theorem}[Br\"and\'en-Krasikov-Shapiro \cite{BKS16}]\label{thm:BKS}
	Let $T:\mathbb{C}[z] \to \mathbb{C}[z]$ be a linear operator defined by
	\begin{equation}\label{eqn:Tp}
		T(p)(z)=q_0(z)p(z)+q_1(z)p(z-1)+\dots+q_k(z)p(z-k),
	\end{equation} where $q_0, \dots, q_k$ are fixed complex-valued polynomials. Then
	$T(\mathcal{HP})\subset\mathcal{HP}$ if and only if $q_i\not\equiv 0$ for at most one $i$, and this $q_i$ is in $\mathcal{HP}$.
\end{theorem}

Recall that the well-known Hermite-Poulain theorem \cite[page 4]{Ob63} states that a finite order linear differential operator $T:=a_0+a_1d/dx+\dots+a_kd^k/dx^k$ preserves the class $\mathcal{HP}$ if and only if its symbol polynomial $Q_T(t)=a_0+a_1t+\dots+a_kt^k$ is hyperbolic. Thus, Theorem \ref{thm:BKS} is a somewhat difference operator analog of the Hermite-Poulain theorem. 
Inspired by this result of Br\"and\'en-Krasikov-Shapiro, we consider finite order linear $q$-difference, difference and differential operators with possibly non-polynomial coefficients,
and obtain three corollaries of Theorem \ref{thm:FD0}.

\begin{corollary}\label{cor:FD}
	Let $q$ be a nonzero real number such that $|q|\neq 1$ and $p\in\mathcal{M}$. Let $T_1:\mathcal{M} \to \mathcal{M}$ be the finite $q$-difference operator defined by $$T_1(p(z)):=a_0p(z)+a_1p(qz)+\dots+a_kp(q^kz),$$ where $a_i\in\mathcal{M}_1$ are entire, for $i=0, 1, \dots, k$. Then, the following statements hold:
	\begin{enumerate}
		\item\label{1} Suppose there exists some $p\in\mathcal{LP}(A;2)$ such that $T_1(p)\in\mathcal{LP}(A;2)$. Then $a_i\not\equiv 0$ for at most one $i$. 
		\item\label{2} The finite $q$-difference operator $T_1$ preserves the class $\mathcal{LP}_0(A; 2)$ if and only if $a_i\not\equiv 0$ for at most one $i$, and this $a_i$ is in the $A$-type I class.
		\item\label{3} Suppose $a_i\in\C[z]$ for all $0\le i \le k$, then $T_1(\mathcal{LP}(\R; 2))\subset\mathcal{LP}(\R; 2)$ if and only if $a_i\not\equiv 0$ for at most one $i$, and this $a_i\in\mathcal{HP}$.
	\end{enumerate}
\end{corollary}

\begin{corollary}\label{cor:FDO1} Let
	$T_2:\mathcal{M} \to \mathcal{M}$ be the finite linear difference operator defined by
	\begin{equation}\label{eqn:T}
		T_2(p(z)):=a_0p(z)+a_1p(z+c_1)+\dots+a_kp(z+c_k)
	\end{equation} where $c_i\in\mathbb{R}\setminus\{0\}$.
	Suppose that each $a_i$ is of the form 
	\begin{equation}\label{eqn:coeff}
		a_i(z)=g_i(z)e^{-\mu_iz^2+d_iz}  \end{equation}
	where $\mu_i$'s are mutually distinct complex numbers with $\mathrm{Re}\,\mu_i>0$,  $d_i\in\mathbb{C}$ and $g_i$ is entire in $\mathcal{M}_1$, for $i=0, \dots, k$. Then the following assertions hold:
	\begin{enumerate}
		\item\label{a} If there exists some $p\in\mathcal{LP}(A)$ such that $T_2(p)\in\mathcal{LP}(A)$, then $a_i\not\equiv0$ for at most one $i$, and this $a_i\in\mathcal{LP}(A; 2)$.
		\item\label{b} Let $t=0$ or $1$, the finite linear difference operator $T_2$ preserves the class $\mathcal{LP}_t(A)$ if and only if $a_i\not\equiv0$ for at most one $i$, and this $a_i\in\mathcal{LP}_t(A; 2)$.
	\end{enumerate} 
\end{corollary}

\begin{corollary}\label{cor:FDO2} Let
	$T_3$ be the finite linear differential operator given by
	\begin{equation}\label{eqn:T2}
		T_3(p(z)):=a_0p(z)+a_1p'(z)+\dots+a_kp^{(k)}(z).
	\end{equation} Suppose that each $a_i$ is of the form (\ref{eqn:coeff}) where $\mu_i$'s are mutually distinct complex numbers with $\mathrm{Re}(\mu_i)>0$,  $d_i\in\mathbb{C}$ and $g_i$ is entire in $\mathcal{M}_1$, for $i=0, \dots, k$. If $T_3(\mathcal{LP}(\R))\subset\mathcal{LP}(\R)$, then $a_i\not\equiv0$ for at most one $i$ and this $a_i\in\mathcal{LP}(\R; 2)$.
\end{corollary}

Corollaries \ref{cor:FD}-\ref{cor:FDO2} give necessary and/or sufficient conditions of such operators preserving some sub-classes of $\mathcal{A}$-entire functions.
In fact, Corollary \ref{cor:FD} could be seen as a generalisation of the Hermite-Poulain theorem (\'a la work \cite{BKS16}) in the $q$-difference context. 
Corollary \ref{cor:FDO1} is a transcendental version of Theorem \ref{thm:BKS} and sort of a complement to the results of \cite{KTV17} and to the Hermite-Poulain theory, developed in \cite{KTV20}, for difference operators of finite order. 
Finally, Corollary \ref{cor:FDO2} is somewhat an extension of the Hermite-Poulain theorem to some sub-classes of $\mathcal{A}$-entire functions.

\subsection{Zeros of differential polynomials}\label{sec:W}
In 1989, Sheil-Small \cite{SSm89} settled a longstanding conjecture of Wiman (1911). 
As a consequence, Sheil-Small obtained that \emph{if $f$ is a real entire function (mapping the real line to itself) of finite order and $ff''$ has no non-real zeros, then $f\in\mathcal{LP}(\R)$.} 
This result also solves a problem posed by Hellerstein (see \cite[p.552, Problem 4.28]{BCl80}). Later, Bergweiler et al.\cite{BEL03} completed  Sheil-Small's result to the real entire function $f$ with infinite order: \emph{for every real entire function $f$ of infinite order, $ff''$ has infinitely many non-real zeros}. Thus, combining the results of Sheil-Small and Bergweiler et al., we have the following
\begin{theorem}[Sheil-Small \cite{SSm89}, Bergweiler-Eremenko-Langley \cite{BEL03}]\label{thm:SSB}
	Let $f$ be a real entire function and $ff''$ has only real zeros, then $f\in\mathcal{LP}(\R)$.
\end{theorem}
Applying a result of Bergweiler et al. to the function of the form $$f=\exp\int_0^zg(t)dt,$$ one can obtain the following result which also follows from a result of Bergweiler and Fuchs \cite{BF93}.
\begin{theorem}[Bergweiler-Fuchs \cite{BF93}; Bergweiler-Eremenko-Langley \cite{BEL03}]\label{cor:BEL}
	For every real transcendental entire function $g$, the function $g'+g^2$ has infinitely many non-real zeros.
\end{theorem}
In 2005, Bergweiler et al. \cite{BEL05} extended Theorem \ref{cor:BEL} to the real meromorphic functions, and considered the zeros of $f'+f^m$ where $m\geq 3$.
It is natural to ask if $f'$ can be replaced by any linear differential polynomial of $f$ or linear difference polynomial of $f$ (Langley \cite[page 108]{langley2009} asked a similar question when $f$ is a real entire function with finitely many non-real zeros). In general, this is not true. For example, let $f=e^{-z}$ which is a real entire function, then $f''+f'+f^m=e^{-mz}$ has no zeros for any integer $m$.  However, applying Theorem \ref{thm:FD0}, we do have a positive result if we restrict to certain sub-classes of real entire functions.

\begin{corollary}\label{cor:SS}
	Let $P$ be a complex polynomial with degree at least two and $P(0)=0$. If $f$ is in class $\mathcal{A}$ with $2\leq\rho(f)<\infty$ and $N(r, 1/f)=o(r)$,
	and $L$ is a non-constant linear differential operator with coefficients in $\mathbb{C}$, then $L(f')+P(f)$ is not in class $\mathcal{A}$. Hence $L(f')+P(f)$ has infinitely many non-real zeros. In particular, this result holds for $f$ in $\mathcal{LP}(A; 2)$ with $N(r, 1/f)=o(r)$.
\end{corollary}

\begin{remark}
	The condition that  $2\leq\rho(f)<\infty$ is necessary as can be seen from the above example for $f=e^{-z}$. It would be interesting to see if the condition $N(r, 1/f)=o(r)$ 
	is also necessary.
\end{remark}
\begin{remark}Essentially the same proof also works when the above differential operator $L(f')$ is replaced by the linear difference operator $a_nf(z+c_n)+\cdots+a_1f(z+c_1)$ where  $a_i\in\C$ and those $c_i$ are mutually distinct nonzero constants.
\end{remark}

\begin{remark}
	In Corollary \ref{cor:SS}, we assume the degree of $P$ is at least two. When $\deg P=1$, Langley \cite[Theorem 1.4]{langley2009} showed that if $f$ is an infinite order real entire function with finitely many non-real zeros, then $f''+ \omega f$ has infinitely many non-real zeros for any positive $\omega$. 
\end{remark}

The rest of the paper is organized as follows. In Sect. \ref{sec:lem} we state several results that will be used in our proofs. Then we prove our main result (Theorem \ref{thm:FD0}) in Sect. \ref{sec:thm}, Corollaries \ref{cor:FD}--\ref{cor:FDO2} in Sect. \ref{sec:cor} and finally Corollary \ref{cor:SS} in Sect. \ref{sec:5}.

\section{Some Lemmata}\label{sec:lem}
As we will apply the half-plane version of Borel's lemma by Rossi \cite{Ro88} (see Lemma \ref{lem:abor}) to prove Theorem \ref{thm:FD0}, we first introduce Tsuji's characteristic of a meromorphic function in the upper (lower) half-plane  (see \cite{Tsu50}). 

Let $\mathfrak{n}_u(t, f)$ be the number of poles of $f$ in $\{z: |z-it/2|\leq t/2, |z|\geq 1\}$, where $f$ is meromorphic in the open upper half plane. Define
$$\mathfrak{N}_u(r, f):=\int_1^r\frac{\mathfrak{n}_u(t, f)}{t^2}dt=\sum_{1\leq r_k\leq r\sin\theta_k}\left(\frac{\sin\theta_k}{r_k}-\frac{1}{r}\right),$$
$$\mathfrak{m}_u(r, f):=\frac{1}{2\pi}\int_{\arcsin(r^{-1})}^{\pi-\arcsin(r^{-1})}\log^+|f(r\sin\theta e^{i\theta})|\frac{d\theta}{r\sin^2\theta},$$ and $$\mathfrak{T}_u(r, f):=\mathfrak{N}_u(r, f)+\mathfrak{m}_u(r, f)$$
where $r_ke^{i\theta_k}$ are the poles of $f$ in $\{z:\mathrm{Im}\, z>0\}$.
Similarly, one can also define $\mathfrak{m}_l(r, f), \mathfrak{N}_l(r, f)$ and $\mathfrak{T}_l(r, f)$ for functions meromorphic in the open lower half plane.
\begin{lemma}[\cite{LO63}]\label{lem:Ros}
	Let $f$ be meromorphic in the open upper (lower) half plane. Define $$m_{\alpha, \beta}(r, f) :=\frac{1}{2\pi}\int_{\alpha}^{\beta}\log^+|f(re^{i\theta})|d\theta.$$ Then $$\int_r^{\infty}\frac{m_{0, \pi}(t, f)}{t^3}dt\leq \int_r^{\infty}\frac{\mathfrak{m}_u(t, f)}{t^2}dt$$
	$$\left(\int_r^{\infty}\frac{m_{\pi, 2\pi}(t, f)}{t^3}dt\leq \int_r^{\infty}\frac{\mathfrak{m}_l(t, f)}{t^2}dt\right).$$
\end{lemma}

\begin{lemma}[\cite{Ro88}]\label{lem:abor}
	Let $n\geq 2$, $G=\{f_0, \dots, f_n\}$ be a set of meromorphic functions in $\mathrm{Im}\, z>0$ such that any proper subset of $G$ is linearly independent over $\mathbb{C}$. If $G$ is linearly dependent over $\mathbb{C}$, then for all positive $r$ except possibly a set of finite measure, $$\mathfrak{T}_u(r) =O\left(\sum_{k=0}^n\left(\mathfrak{N}_u(r, f_k)+\mathfrak{N}_u\left(r, \frac{1}{f_k}\right)\right)+\log \mathfrak{T}_u(r)+\log r\right),$$
	where $\mathfrak{T}_u(r):=\max\{\mathfrak{T}_u(r, f_i/f_j)|\, 0\leq i, j\leq n\}$.
\end{lemma}

We also need the following generalisation \cite{nevanlinna1929} of Borel's lemma \cite{Bor21}.  
\begin{lemma}[Corollary 4.5 of \cite{BCL93}]\label{lem:BCL}
	Let $a_i$ and $g_i$, $i=1, \dots, k$ be nonzero meromorphic and entire functions in $\C$ respectively, satisfying $$a_1e^{g_1}+a_2e^{g_2}+\dots+a_ke^{g_k}\equiv 0.$$ If $T(r, a_j)=o(T(r, e^{g_m-g_l}))$, for any $m\neq l$ and $1\leq j\leq k$, then $a_j(z)\equiv 0$ for all $1\leq j\leq k$.
\end{lemma}

\section{Proof of  Theorem \ref{thm:FD0}}\label{sec:thm}

We may assume that $a_i \neq 0$ for any $i=1,\dots,n$ (otherwise,  we can relabel $a_i$ so that we can replace $n$ by a smaller number). Let $g_0=F$ and $g_i=a_if_i$ for $i=1, \dots, n$.
By the assumption that $f_1, \dots, f_n$ are linearly independent over $\mathcal{M}_2$ and $a_i \in \mathcal{M}_1 \subset \mathcal{M}_2$, we have $g_1, \dots, g_n$ are linearly independent over $\mathbb{C}$.
Consider the set $G=\{g_0, \dots, g_n\}$, then $G$ is linearly dependent over $\mathbb{C}$ and any proper subset of $G$ is linearly independent over $\mathbb{C}$. Therefore, $G$ satisfies the assumptions of Lemma \ref{lem:abor}.

Suppose that $F=a_1f_1+\dots+a_nf_n \in \mathcal{A}\backslash\{0\}$. Since all $g_0, f_1, \dots, f_n$ are in $\mathcal{A}$ (hence entire), it is easy to check from the definitions of class $\mathcal{A}$ and $\mathfrak{N}_u(r, *)$ that $$\mathfrak{N}_u(r, 1/g_0)=O(1)\ \mathrm{and}\ \mathfrak{N}_u(r, 1/f_i)=O(1).$$ We
also have for $i=1, \dots, n$, $\mathfrak{N}_u(r, g_i)=0$ because each $a_i$ is entire.

Recall that $T(r, a_i)=o(r)$. Then for each $i=1, \dots, n$,
\begin{align*}
	\mathfrak{N}_u(r, 1/g_i)+\mathfrak{N}_u(r, g_i)
	&\le \mathfrak{N}_u(r, 1/a_i)+ \mathfrak{N}_u(r, 1/f_i) +O(1)\\
	&\le N(r,1/a_i)+O(1) \le T(r,a_i)=O(r^{\epsilon}) 
\end{align*}
for some positive $\epsilon<1$. We also have $\mathfrak{N}_u(r, 1/g_0)+\mathfrak{N}_u(r, g_0) =O(1)$. Therefore, we can deduce from Lemma \ref{lem:abor}, that $$\mathfrak{T}_u(r)=O(r^{\epsilon}),\ \text{and hence}\ \mathfrak{T}_u(r, g_i/g_j)=O(r^{\epsilon})$$ for $i, j=0, \dots, n$.
From the definition of $\mathfrak{T}_u(r, g_{i}/g_{j})$, it follows that $$\mathfrak{m}_u(r, g_{i}/g_{j})=O(r^{\epsilon}).$$
Similarly,  $$\mathfrak{m}_l(r, g_{i}/g_{j})=O(r^{\epsilon}).$$

Notice that $m(t, g_{i}/g_{j})=m_{0, \pi}(t, g_i/g_j)+m_{\pi, 2\pi}(t, g_i/g_j)$. Then Lemma \ref{lem:Ros} implies that for any $r>0$, 
\begin{align*}
	\frac{m(r, g_{i}/g_{j})}{2r^2}&\leq \int_{r}^{\infty}\frac{m(t, g_{i}/g_{j})}{t^3}dt\leq O(r^{\epsilon-1})
\end{align*}
and hence $m(r, g_i/g_j)=o(r^2)$.

Using the fact that $a_i \in \mathcal{M}_1, N(r,f_i)=o(r^2)$ and $a_i,f_i$ are entire, we have $$N(r, g_i/g_j)\leq N(r, 1/a_j) + N(r, 1/f_j) \le T(r, 1/a_j) + N(r, 1/f_j) =o(r^2).$$ Therefore, $T(r, g_i/g_j)=o(r^2)$.

Take $i=0$ and $j=1$, then $T(r, \frac{F}{a_1f_1})=o(r^2)$. Let $b=\frac{F}{a_1f_1}$ so that $F=ba_1f_1$ and $T(r,(1-b)a_1)=o(r^2)$. Hence
$$(1-b)a_1f_1+a_2f_2+\dots+a_nf_n=0.$$ As $f_1, \dots, f_n$ are linearly independent over $\mathcal{M}_2$, we must have $a_i=0$ for each $i \neq 1$ and $(1-b)a_1=0$. Thus, the only possibility is $n=1$ and $b=1$ so that $F=a_1f_1$ and $a_1\in\mathcal{A}$.

The converse assertion is clearly true and this completes the proof.

\section{Proof of  Corollaries \ref{cor:FD},  \ref{cor:FDO1} and  \ref{cor:FDO2}}\label{sec:cor}
\subsection{Proof of Corollary \ref{cor:FD}$(\ref{1})$}
Since $p\in\mathcal{LP}(A;2)$, we have $p(z)=h(z)e^{-\alpha z^2+\beta z}$ where $\alpha >0$ and $h(z)$ is of the form (\ref{eqn:f})  so that $h\in \mathcal{A}\backslash\{0\}$ and $T(r, h)=o(r^2)$.
Let $b_i(z)=h(q^iz)$ and $f_i(z)=b_i(z)e^{-\alpha q^{2i}z^2+\beta q^iz}$ for $i=0, \dots, k$ where $|q|\neq 0, 1$.

Clearly $f_i\in\mathcal{A}\backslash\{0\}$ if $b_i\in\mathcal{A}\backslash\{0\}$. To see $b_i\in\mathcal{A}\backslash\{0\}$, let $\{\beta_k\}$ and $\{z_k\}$ be zeros of $b_i$ and $h$ respectively. Since $q$ is real, it follows that
\begin{align*}
	\sum|\mathrm{Im} \frac{1}{\beta_k}|=|q|^i\sum|\mathrm{Im} \frac{1}{z_k}|<\infty
\end{align*} and hence $b_i \in \mathcal{A}$. Since $\mathcal{LP}(A;2) \subset \mathcal{A}$, each $f_i$ and $T_1(p)=a_0f_0+\cdots+a_kf_k$ are in class $\mathcal{A}\backslash\{0\}$. Finally, we notice that $N(r,1/f_i)=N(r,1/b_i)=N(|q|^{i}r,1/h)+O(1)\le T(|q|^{i}r,h)+O(1)=o(r^2)$.

To show that at most one $a_i \not\equiv 0$ by applying Theorem \ref{thm:FD0}, 
it remains to prove that $f_0, \dots, f_k$ are linearly independent over $\mathcal{M}_2$.
Suppose there exist $c_i \in \mathcal{M}_2$, for $i=0, \dots, k$, 
such that $c_0f_0+\cdots +c_kf_k=0$. Let $g_i=e^{-\alpha q^{2i}z^2+\beta q^iz}$. As $$\frac{g_i}{g_j}=e^{-(\alpha q^{2i}-\alpha q^{2j})z^2+(\beta q^i-\beta q^j)z}$$ and $\alpha q^{2i}-\alpha q^{2j}\neq 0$ for $i\neq j$, we have whenever $i\neq j$, $$T(r, g_i/g_j)\geq Cr^2,\ \text{as}\ r\rightarrow\infty,\ C>0.$$ 
On the other hand, one can check easily that $T(r,h(q^jz))=T(|q|^jr, h)+O(1)$ (see \cite[Page 249]{bergweiler1998}) so that $T(r,h(q^jz))=o(r^2)$ for all $0\le j \le k$ and hence $T(r, c_jb_j)=o(r^2)$ as $c_j \in \mathcal{M}_2$. Since 
$$c_0b_0e^{-\alpha z^2+\beta z}+\cdots +c_kb_ke^{-\alpha q^{2k}z^2+\beta q^kz}=c_0f_0+\cdots +c_kf_k=0,$$
by Lemma \ref{lem:BCL}, $c_jb_j=0$ for all $j=0, \dots, k$ which implies that $c_j=0$ for all $j$ and therefore $f_0, \dots, f_k$ are linearly independent over $\mathcal{M}_2$.

\subsection{Proof of Corollary \ref{cor:FD}$(\ref{2})$}
Since $T_1(\mathcal{LP}_0(A;2))\subset\mathcal{LP}_0(A;2)\subset\mathcal{LP}(A;2)$, it follows from part one that $a_i\not\equiv 0$ for at most one $i$.
Therefore, it remains to prove that this $a_i$ is in class $A$-type I. 

Without loss of generality, we may assume that $i=0$ and hence $T_1(p(z))=a_0(z)p(z)=a_0(z)h(z)e^{-\alpha z^2+\beta z}$, where $h(z)=cz^n\displaystyle\prod_{k=1}^{\infty}(1-z/z_k) \in \mathcal{M}_1$, $c\in \mathbb{R}\backslash\{0\}$. 
Since $T_1(p(z))\in\mathcal{LP}_0(A, 2)$, we have \begin{equation}\label{eqn:lp0}
	a_0(z)h(z)e^{-\alpha z^2+\beta z}=h_1(z)e^{-\alpha_1 z^2+\beta_1 z},
\end{equation} where $h_1(z)=c_1z^m\displaystyle\prod_{k=1}^{\infty}(1-z/w_k) \in \mathcal{M}_1$, $c_1\in \mathbb{R}\backslash\{0\}$. Since $a_0h$ and $h_1$ are in $\mathcal{M}_1$, we can apply Lemma \ref{lem:BCL} to conclude that $a_0h=h_1$ and hence 
$$a_0=c'z^l\displaystyle\prod_{w_k\neq z_k}(1-z/w_k)\quad\mathrm{for\ some}\ c'\in\R\setminus\{0\}\ \mathrm{and}\ l\in\N.$$ 
Since $$\sum_{w_k\neq z_k}|w_k|^{-1}\leq\sum_{k=1}^{\infty}|w_k|^{-1}<\infty \quad \mathrm{and} \quad  \sum_{w_k\neq z_k}|\mathrm{Im}\frac{1}{w_k}|\leq \sum_{k=1}^{\infty}|\mathrm{Im}\frac{1}{w_k}|<\infty,$$ it follows that $a_0$ is in class $A$-type I. 

The converse assertion is obvious and we complete the proof.
\subsection{Proof of Corollary \ref{cor:FD}$(\ref{3})$}
Since $a_i\in\C[z] \subset \mathcal{M}_1$ for all $i$, and $T_1(\mathcal{LP}(\R; 2))\subset\mathcal{LP}(\R; 2)$, by Corollary \ref{cor:FD}(\ref{1}), we conclude that $a_i\not\equiv 0$ for at most one $i$ and there is no harm to assume that only $a_0\not\equiv 0$. The property that $a_0\in\mathcal{HP}$ then follows from the identity (\ref{eqn:lp0}).

\subsection{Proof of Corollary  \ref{cor:FDO1}$(\ref{a})$}
Let $p=h(z)e^{-\alpha z^2+\beta z} \in\mathcal{LP}(A)$ be expressed in the form of (\ref{eqn:f}) where $T(r,h)=o(r^2)$. Let $c_0=0$ so that we can write $T_2(p)$ as 
$$T_2(p)=b_0f_0+b_1f_1+\dots+b_kf_k$$ where $b_i=g_i$ and $f_i(z)=h(z+c_i)e^{-(\alpha (z+c_i)^2+\mu_iz^2)+\beta (z+c_i)+d_iz},$ with $\alpha+\mu_i\neq 0$ for all $i$ because $\mathrm{Re}(\mu_i)>0$ and $\alpha\geq0$.  
Since $N(r,1/h)\le T(r,h) +O(1)=o(r^2)$, $\lambda(h) \le 2$. By \cite[Theorem 2.2]{CF08}, it follows that $N(r,1/f_i)=N(r,1/h(z+c_i))=N(r,1/h)+O(r^{\lambda(h)-1 + \epsilon}) + O(\log r)$ for any positive $\epsilon$. Hence $N(r,1/f_i))=o(r^2)$ for all $i$. 

Since each $c_i$ is real, one can check that if $|z_k|>2|c_i|$, then $$|\mathrm{Im} \frac{1}{z_k-c_i}|\le 4|\mathrm{Im}\, \frac{1}{z_k}|$$ and hence $f_i\in\mathcal{A}$. Now suppose $T_2(p) \in \mathcal{LP}(A)$. In order to apply Theorem \ref{thm:FD0} to  show that at most one $a_i \not\equiv 0$, we only need to show that $f_0, \dots, f_k$ are linearly independent over $\mathcal{M}_2$. This suffices to show that each $h(z+c_i)$ is in $\mathcal{M}_2$ because by  
Lemma \ref{lem:BCL} and the fact that $\mu_i$'s are distinct, $\{e^{-(\alpha (z+c_i)^2+\mu_iz^2)+\beta (z+c_i)+d_iz}: i=0,\ldots,k\}$ is linearly independent over $\mathcal{M}_2$.

Since $T(r,h)=o(r^2)$, its order $\sigma$ is at most two. By \cite[Theorem 2.1]{CF08}, it follows that $T(r,h(z+c_i))=T(r,h)+O(r^{\sigma-1 + \epsilon}) + O(\log r)$ for any positive $\epsilon$. Hence $T(r,h(z+c_i))=o(r^2)$ for all $i$ and we are done. 

\subsection{Proof of Corollary  \ref{cor:FDO1}$(\ref{b})$}
The argument is similar to that of Corollary \ref{cor:FD}(\ref{b}),  and we omit the details of the proof.

\subsection{Proof of Corollary \ref{cor:FDO2}} 
Let $p(z)=h(z)e^{-\alpha z^2+\beta z}\in\mathcal{LP}(\R)$, where $h$ is given by (\ref{eqn:f}) with $T(r,h)=o(r^2)$. Then $p^{(i)}(z)=h_i(z)e^{-\alpha z^2+\beta z}$, where $h_i(z)=h_{i-1}'(z)+(-2\alpha z+\beta )h_{i-1}(z)$ and $h_0(z)=h(z)$. By  \cite[Theorem A]{HW77}, it follows that $p^{(i)}\in \mathcal{LP}(\R)\subset \mathcal{A}$. Hence $f_i:=p^{(i)}e^{-\mu_iz^2+d_iz}=h_ie^{-(\alpha + \mu_i)z^2+(\beta +d_i)z} \in \mathcal{A}$. As $T(r,h_0)=o(r^2)$ and $T(r, h_i)=O(T(r, h_{i-1}))$, we have $T(r,h_i)=o(r^2)$ for all $i$. Therefore, $N(r,1/f_i)=o(r^2)$ for all $i$.
The rest can be done as in the proof of Corollary \ref{cor:FDO1}(\ref{a}).

\section{Proof of Corollary \ref{cor:SS}}\label{sec:5}
If $L(f')+P(f)\equiv 0$, then $L(f')+P(f)\notin \mathcal{A}$ and we are done. So we may assume that $L(f')+P(f)\not\equiv 0$ and write it in the following form
\begin{align*}
	L(f')+P(f)&=\sum_{k=1}^na_kf^{(k)}+\sum_{m=1}^lb_mf^m=(\sum_{k=1}^na_k\frac{f^{(k)}}{f})f+\sum_{m=1}^lb_mf^m\\
	&=\sum_{k=1}^lc_kf^k
\end{align*}
where $a_i \in \mathbb{C}$, $c_1=b_1+\sum_{k=1}^na_k\frac{f^{(k)}}{f}$, $c_m=b_m  \in \mathbb{C}$ for $m=2, \dots, l$. Let $f_j=f^j$ for $j=1, \dots, l$, then $$f_j\in\mathcal{A} \quad  \mathrm{and}\quad N(r,1/f_j)=jN(r, 1/f)=o(r) $$
and hence $N(r,1/f_j)=o(r^2)$. As $\rho(f)\geq 2$, by Valiron's Theorem (\cite[Theorem 2.2.5]{IIpo11}), $f_1, \dots, f_l$ are linearly independent over $\mathcal{M}_2$. Since $f$ is entire with $\rho(f)<\infty$ and $N(r, 1/f)=o(r)$,
it follows from the logarithmic derivative lemma and $N(r, f^{(k)}/f)\le N(r,1/f)$ that
\begin{align*}
	T(r, f^{(k)}/f)&=m(r, f^{(k)}/f)+N(r, f^{(k)}/f)\\
	&= o(r).
\end{align*} 
Hence  $f^{(k)}/f\in \mathcal{M}_1 $ and so is $c_1$.
Now suppose $L(f')+P(f)=c_1f+\cdots+c_kf^k \in \mathcal{A}\backslash\{0\}$. By Theorem \ref{thm:FD0}, we have $c_i\not\equiv0$ for at most one $i$.  As $\deg P\geq 2$,  $b_l\neq 0$,
therefore, $c_1\equiv 0$, i.e., $L(f')+b_1f=a_nf^{(n)}+\dots+a_1f'+b_1f=0$. This implies that $\rho(f)\leq 1$ which is a contradiction to the order of growth of $f$ is at least $2$.

\section*{Acknowledgments}
The authors would like to thank the anonymous
referees for their valuable suggestions and helpful comments. 

  \bibliographystyle{abbrv}
\bibliography{mybibfile}
\end{document}